\documentclass[12pt,draft,amsmath,amscd,amsbsy,amssymb]{amsart}
 \relax

\chardef\bslash=`\\ 

\makeatletter \def\verbatim{\interlinepenalty\@M \@verbatim
\leftskip\@totalleftmargin\advance\leftskip2pc \frenchspacing\@vobeyspaces
\@xverbatim} \makeatother \hfuzz1pc

\makeatletter \def\dgt@k{\dg@DX=-3 \dg@DY=2 \dg@SIZE=3} \makeatother

\makeatletter \def\dgt@kk{\dg@DX=3 \dg@DY=-1 \dg@SIZE=3}

\theoremstyle{plain}
\newtheorem{thm}{Theorem}[section]

\newtheorem{lemma}[thm]{Lemma}

\theoremstyle{definition}

\numberwithin{equation}{section}

\newcommand{\comp}{{\mathrm C}{\mathrm o}{\mathrm m}{\mathrm p} }

\usepackage[all]{xy}

\begin{document}


\title[]
{Open-multicommutativity of the probability measure functor}
\author{R. Kozhan}

\address{} \email{}

\author{M.~Zarichnyi}

\address{Department of Mechanics and Mathematics, Lviv National University,
Universytetska 1, 79000 Lviv, Ukraine} \email{topology@franko.lviv.ua,
mzar@litech.lviv.ua}
\thanks{}

\begin{abstract} For the functors acting in the category of compact Hausdorff
spaces, we introduce the so-called open multi-commutativity property, which
generalizes both bicommutativity and openness, and prove that this property is
satisfied by the functor of probability measures.
\end{abstract}

\maketitle

\section{Introduction}\label{s:intro}

It is well-known that the construction of space of probability measures $P$ is
functorial in the category $\comp$ of compact Hausdorff spaces. The functor $P$
is normal in the sense of E.V. Shchepin \cite{Sh}. It is well-known that the
functor $P$ is open, i.e. it preserves the class of open surjective maps. This
was first proved by Ditor and Eifler \cite{DE}. E.V. Shchepin \cite{Sh}
discovered tight relations between the properties of openness and
bicommutativity. In particular, he proved that every open functor in $\comp$ is
bicommutative, i.e. preserves the class of the bicommutative diagrams in the
sense of Kuratowski.

In this paper we introduce the so-called open multi-commutativity property and
show that this property is satisfied by the functor $P$.

\section{Preliminaries}

\subsection{Openness} We say that a functor in $\comp$ is {\it open} if it
preserves the class of open surjective maps.

\subsection{Bicommutativity} A commutative diagram
$$\begin{array}{cc}
{\mathcal D}=& {\xymatrix{X\ar[r]^f\ar[d]_g& Y\ar[d]\\ Z\ar[r]& T}}
\end{array}$$ in $\comp$ is
said to be {\it bicommutative} if its characteristic map $\chi_{\mathcal
D}=(f,g)\colon X\to Y\times_TZ$ is onto. We say that a functor in $\comp$ is {\it
bicommutative} if it preserves the bicommutative diagrams (see \cite{Sh}).

\subsection{Open-multicommutativity}

Suppose that $O$ is a finite partially ordered set. We regard $O$ as a category
and consider a functor $F\colon O\to\comp$ (i.e. a diagram in $\comp$ indexed by
the objects of $O$). A {\it cone} over $F$ consists of a space $X\in|\comp|$ and
a family of maps $\{X\to Fd\}_{d\in |O|}$ that satisfy obvious commutativity
conditions. Given such a cone, $\mathcal C=(\{X\to Fd\}_{d\in |O|})$, we denote
by $\chi_{\mathcal C}\colon X\to\lim F$ its characteristic map.

We say that the cone $\mathcal C$ is {\it open-multicommutative} if its
characteristic map is an open onto map.

We say that a functor in $\comp$ is {\it open-multicommutative} if it preserves
the class of open-multicommutative diagrams.

Note that, if $|O|$ consists of one object then the open-multicommutativity
reduces to openness.

 \section{Main result}

In the sequel, we need a more detailed description of the diagrams under
consideration. We denote the spaces in the diagram by $X_i$, where $i$ is a
generic element of a partially ordered set. If $i\ge j$, then we denote by
$\varphi_{ij}$ the map from $X_i$ to $X_j$. The limit, $\lim D$, of the diagram
$D$ can be naturally embedded into the set $$\prod\{ X_i\mid i\text{ is a maximal
element of the set of indices}\}.$$ We will denote this embedding by $h$.

One can naturally define the category of diagrams in $\comp$ with the same index
set $O$. We denote this category by $\comp^O$.

\begin{lemma}\label{} The operations of the limits of the inverse systems
and the limits of the diagrams commute.
\end{lemma}

\begin{proof} Straightforward.
\end{proof}

Let $D$ be a diagram in $\comp^O$. Without loss of generality, we may suppose
that the set of maximal elements of $O$ is $\{1,\dots,k\}$. Then $\lim D\subset
X_1\times\dots\times X_k$.

Consider the diagram $P\left( D\right) .$ The limit of this diagram is $${\lim
}P\left( D\right) \subset P\left( X_1\right)\times...\times P\left( X_k\right)
.$$ There exists the unique map $\chi \colon P(\lim D)\to\lim P(D)$.

\begin{thm}
The map $\chi $ is open.
\end{thm}

\begin{proof}

 We will need the following lemmas.
\begin{lemma}\label{l:2} Let maps $f_{X_i}:X_i\rightarrow X_i^{\prime }$ be such
that the diagrams
\begin{equation}
\xymatrix{X_i\ar[r]^{\varphi _{ij}}\ar[d]_{f_{X_i}}&X_j\ar[d]^{f_{X_j}}\\
X_i^{\prime }\ar[r]_{\varphi _{ij}^{\prime }}&X_j^{\prime }}
\end{equation}are commutative for every $i\ge j$. In other words, $(f_{X_i})$ is
a morphism in $\comp^O$ of the diagram $D$ to a diagram $D'$ in which the spaces
and maps are endowed with ``prime''. Then the diagram
\[
\xymatrix{P\left( {\lim }D \right)\ar[r]^\chi\ar[d]_{P\left( h\circ
\overset{k}{\underset{i=1}{\prod }}f_{X_i}\circ h^{-1}\right)
}&{\lim } P\left(D\right)\ar[d]^{h\circ \overset{k}{\underset{i=1}{\prod
}}P\left( f_{X_i}\right)
\circ h^{-1}}\\ P\left( {\lim }D '\right)\ar[r]_{\chi ^{\prime }}&
{\lim } P\left(D'\right)}
\]
is bicommutative.
\end{lemma}

\begin{proof}
In order to prove the bicommutativity of this diagram we should, given $\tau
_0\in P\left(D\right) $ and $\mu =\left(
\mu _1,...,\mu _k\right) \in {\lim } P\left(D
\right)  $ such that
\[\chi ^{\prime }\left( \tau _0\right)
=\overset{k}{\underset{i=1}{\prod }}P\left( f_{X_i}\right) \left(
\mu \right) =\left( Pf_{X_1}\left( \mu _1\right)
,...,Pf_{X_i}\left( \mu _k\right) \right), \] find a measure $\tau
\in P\left(D \right)
$ for which
$
\chi \left( \tau \right) =\mu $ and $P\left(
\overset{k}{\underset{i=1}{\prod }}f_{X_i}\right) \left( \tau
\right) =\tau _0.
$

Let us define diagrams $D_i,$ $i\in \left\{ 1,...,k\right\} $ in the following
way. We replace every $X_j$, $j=i+1,\dots,k$, by $X'_j$ and every $\varphi_{js}$
by $\varphi_{js}f_{X_j}$. Every other maps and spaces are determined by the
following conditions.

1. If $ s\ge t $ and $X_s^{\prime }\in D_i$ then $X_t^{\prime }\in D_i.$

2. Only one of $X_j$ or $X_j^{\prime }$ can be in $D_i$ for every $j.$

3. If $s\ge t$ and $X_s\in D_i$ then $X_t\in D_i.$

4. For $s\ge t$ if $X_s,X_t\in D_i$ then $\varphi _{st}\in D_i,$ if $X_s^{\prime
},X_t^{\prime }\in D_i$ then $\varphi _{st}^{\prime }\in D_i,$ if
$X_s,X_t^{\prime }\in D_i$ then $\varphi _{st}^{\prime }\circ f_{X_t}\in D_i$.

Now for every $i\in \left\{ 1,...,k\right\} $ denote by ${\mathcal D}^i$ the
square diagram
\[
\xymatrix{{\lim }{D_i}\ar[rrrr]^{\pi _i}  \ar[d]_{H^i=h\circ
\left( \underset{j=1}{\overset{i-1}{\prod }}1_{X_j}\times
f_{X_i}\times \overset{k}{\underset{j=i+1}{\prod }}1_{X_j^{\prime }}\right) \circ
h^{-1}}&&&& X_i\ar[d]^{f_{X_i}}\\ {\lim }{D_{i-1}}\ar[rrrr]_{\pi
_i}&&&&{X_i}^{\prime }}
\]This diagram is well-defined. Indeed, let $\left( x_j\right)
\in {\lim }{D_i} $ and
\begin{eqnarray*}
H\left( \left( x_j\right)_j \right) &=&h\left( \left(
\overset{i-1}{\underset{j=1}{\prod }}1_{X_J}\times f_{X_i}\times
\overset{k}{\underset{j=i+1}{\prod }}\right) \left( h^{-1}\left( x_j\right)_j \right)
\right)\\&=&\left( y_j\right) _{j }\notin {\lim }{D_{i-1}}
.
\end{eqnarray*}
For every $s\notin \left\{1,\dots, k\right\} $ we consider arbitrary index $t\ge
s.$ Without loss of generality we can assume that $t\in \left\{1,\dots,
k\right\}$. There are five possibilities.

1. $X_t,X_s\in D_i$ and $X_t,X_s\in D_{i-1}.$ Then $H_t^i\left(
x_t\right) =y_t=x_t\in X_t$ and $H_s^i\left( x_s\right)
=y_s=x_s\in X_s.$ But $\varphi _{ts}\left( y_t\right) =\varphi
_{ts}\left( x_t\right) =x_s=y_s$ therefore $\left( y_j\right)
_{j }\in {\lim
}{D_{i-1}}$.

2. $X_t^{\prime },X_s^{\prime }\in D_i$ and $X_t^{\prime
},X_s^{\prime }\in D_{i-1}.$ Then $H_t^i\left( x_t\right)
=y_t=x_t\in X_t^{\prime }$ and $H_s^i\left( x_s\right) =y_s=x_s\in
X_s^{\prime }.$ But $\varphi _{ts}\left( y_t\right) =\varphi
_{ts}^{\prime }\left( x_t\right) =x_s=y_s$ therefore $\left( y_j\right)
_{j }\in {\lim
}{D_{i-1}}$.

3. $X_t,X_s^{\prime }\in D_i$ and $X_t,X_s^{\prime }\in D_{i-1}.$
Then $H_t^i\left( x_t\right) =y_t=x_t\in X_t^{\prime }$ and
$H_s^i\left( x_s\right) =y_s=x_s\in X_s^{\prime }.$ But
\[\left( \varphi _{ts}\circ f_{X_t}\right)
\left( y_t\right) =\left( \varphi _{ts}\circ f_{X_t}\right) \left(
x_t\right) =x_s=y_s \]
therefore $\left( y_j\right)
_{j }\in {\lim
}{D_{i-1}}$.

 4. $X_t,X_s\in D_i$ and $X_t^{\prime },X_s^{\prime }\in D_{i-1}.$
Then $H_t^i\left( x_t\right) =y_t=f_{X_t}\left( x_t\right) \in X_t^{\prime }$ and
$H_s^i\left( x_s\right) =y_s=f_{X_s}\left( x_s\right) \in X_s^{\prime }.$ Since
diagrams (1) are commutative
\[
\varphi _{ts}^{\prime }\left( y_t\right) =\varphi _{ts}^{\prime
}\left( f_{X_t}\left( x_t\right) \right) =f_{X_s}\left( \varphi
_{ts}\left( x_t\right) \right) =f_{X_s}\left( x_s\right) =y_s
\]therefore $\left( y_j\right)
_{j }\in {\lim
}{D_{i-1}}$.

5. $X_t,X_s\in D_i$ and $X_t,X_s^{\prime }\in D_{i-1}.$ Then $H_t^i\left(
x_t\right) =y_t=x_t\in X_t$ and $H_s^i\left( x_s\right) =y_s=f_{X_s}\left(
x_s\right) \in X_s^{\prime }$ and $\varphi _{ts}^{\prime }\circ f_{X_t}\in
D_{i-1}.$ The commutativity of diagrams (1) implies
\[
\left( \varphi _{ts}^{\prime }\circ f_{X_t}\right) \left(
y_t\right) =f_{X_s}\left( \varphi _{ts}\left( x_t\right) \right)
=f_{X_s}\left( x_s\right) =y_s
\]therefore $\left( y_j\right)
_{j }\in {\lim
}{D_{i-1}}$.

These are all possibilities which can happen and this implies that $\left(
y_j\right) _{j}\in {\lim }{D_{i-1}}$.
 Thus, the diagram is
well defined.

Every diagram $\mathcal D^i$ is bicommutative. Since $P\pi _1\left( \tau
_0\right) =Pf_{X_1}\left( \mu _1\right) ,$ applying the functor
$P$ to the diagram $D^1 $ one can find $\tau _1\in P\left(\lim D_1\right) $ such
that
\[
P\pi
_1\left( \tau _1\right) =\mu _1,\text{ }PH^1\left( \tau _1\right)
=P\left( h\circ \left( f_{X_1}\times
\overset{k}{\underset{j=2}{\prod }}1_{X_j^{\prime }}\right) \circ
h^{-1}\right) \left( \tau _1\right) =\tau _0.
\]We assume that for every $i\in \left\{ 1,...,k\right\} $ we can
define $\tau _i\in P\left( {\lim } {D_i}\right)$ such that
\[
P\pi _i\left( \tau _i\right)
=\mu _i,\text{ }PH^i\left( \tau _i\right) =P\left( h\circ \left(
\overset{i-1}{\underset{j=1}{\prod }}1_{X_J}\times f_{X_i}\times
\overset{k}{\underset{j=i+1}{\prod }}\right) \circ h^{-1}\right)
\left( \tau _1\right) =\tau _{i-1}.
\]This holds due to the fact that
\begin{eqnarray*}
Pf_{X_i}\left( \mu _i\right) &=&P\pi _i\left( \tau _0\right) =P\pi
_i\left( H^1\left( \tau _1\right) \right) =P\pi _i\left( H^1\left(
H^2\left( \tau _2\right) \right) \right)\\&=&P\pi _i\left(
H^1\circ H^2\circ ...\circ H^{i-1}\right) \left( \tau
_{i-1}\right) =P\pi _i\left( \tau _{i-1}\right) .
\end{eqnarray*}
Consider now the map $P\left( h\circ
\overset{k}{\underset{j=1}{\prod }}f_{X_j}\circ h^{-1}\right) .$
Since
\begin{eqnarray*}
\overset{k}{\underset{j=1}{\prod }}f_{X_j} =\left( f_{X_1}\times
\overset{k}{\underset{j=2}{\prod }}1_{X_j}\right) &\circ& \left(
1_{X_1}\times f_{X_2}\times \overset{k}{\underset{j=3}{\prod
}}1_{X_j^{\prime }}\right) \circ ...\\&\circ& \left(
\overset{s-1}{\underset{j=1}{\prod }}1_{X_J}\times f_{X_s}\times
\overset{k}{\underset{j=s+1}{\prod }}1_{X_j^{\prime }}\right)
\circ ...\circ \left( \overset{k-1}{\underset{j=1}{\prod
}}1_{X_J}\times f_{X_k}\right),
\end{eqnarray*} we have
\begin{eqnarray*}
P\left( h\circ \overset{k}{\underset{j=1}{\prod }}f_{X_j}\circ
h^{-1}\right) \left( \tau _k\right)& =&P\left( H^1\circ H^2\circ
...\circ H^k\right) \left( \tau _k\right)\\&=&P\left( H^1\circ
H^2\circ ...\circ H^{k-1}\right) \left( \tau _{k-1}\right)
=...=P\left( H^1\right) \left( \tau _1\right) =\tau _0.
\end{eqnarray*}As it has been proved before
\[
\chi \left( \tau
_k\right) =\left( P\pi _1,...,P\pi _k\right) \left( \tau _k\right)
=\left( \mu _1,...,\mu _k\right) .
\]The measure $\tau =\tau _1$ is the measure we were looking for.

\end{proof}

\begin{lemma}
The map $\chi $ is surjective.
\end{lemma}

\begin{proof}

Denote by
\[
\Gamma _1=\left\{ j\notin \left\{1,\dots,k\right\} \mid 1\ge
j\right\}
\]
\[
\Gamma _2=\left\{ j\notin \left\{1,\dots,k\right\} \mid 2\ge
j\right\} \setminus \Gamma _1
\]
\[
\vdots
\]
\[
\Gamma _k=\left\{ j\notin \left\{1,\dots,k\right\} \mid k\ge
j\right\} \setminus \left( \Gamma _1\cup ...\cup \Gamma
_{k-1}\right) .
\]We denote the indices in $\Gamma _l$ as
$j_1^l,...,j_{m_1}^l$, where $l\in \left\{ 1,...,k\right\} .$ Let $\left( \mu
_1,...,\mu _n\right) \in {\lim }P({D}) .$

We have to find a measure $\tau \in P\left( {\lim }{D}
\right) $ such that $\chi \left(
\tau \right)
=\left( \mu _1,...,\mu _n\right) .$ Denote by $\mathcal D_{11}$ the
following bicommutative diagram
\[
\xymatrix{X_1\times X_{j_1^1}\ar[r]^{\pi _1}\ar[d]_{\pi _1}
&X_1\ar[d]^{\varphi _{1...j_1^1}}\\
X_1\ar[r]_{\varphi _{1...j_1^1}}&X_{j_1^1}}
\]Since $\left( \mu _1,...,\mu _n\right) \in \underset{D}{\lim
}\left\{ P\left( X_1\right) ,...,P\left( X_k\right) \right\} $, applying the
functor $P$ to this diagram we can find a measure $\tau _{11}$ for which $P\pi
_1\left( \tau _{11}\right)
=\mu _1$ and $P\pi _2\left( \tau _{11}\right) =\mu _{j_1^1}.$ Next
we define the diagram $D_{it},$ where $i\in \left\{
1,...,k\right\} $ and $j_t^i\in \Gamma _i$
as
\[
\xymatrix{X_1\times X_{j_1^1}\times ..\times X_{j_1^{m_1}}\times
..\times X_i\times ..\times X_{j_t^i}\ar[rrrrr]^{\pi
_i}\ar[d]_{\pi
_{12...j_{t-1}^i}} &&&&&X_i\ar[d]^{\varphi _{i...j_t^i}}\\
 X_1\times X_{j_1^1}\times
..\times X_{j_1^{m_1}}\times ..\times X_i\times ..\times
X_{j_{t-1}^i}\ar[rrrrr]_{\varphi _{i..j_t^i}\circ \pi
_{j_{t-1}^i}}&&&&& X_{j_t^i}}
\]Applying in natural order for all these diagrams $D_{it},$ where
$i$ runs from $1$ to $k$ and $j_t^i\in \Gamma _i$ we find at the very end a
measure $\tau _n=\tau \in P\left({\lim }{D}
\right) .$ Then for every $i\in
\left\{ 1,...,k\right\} $ we have $P\pi _i\left( \tau \right)
=P\pi _i\circ P\pi _{12...i}\circ P\pi _{12...i+1}\circ ...\circ
P\pi _{12...n-1}\left( \tau \right) =\mu _i.$ For each $j\notin
\left\{1,\dots, k\right\} $ there exists $i\in \left\{
1,...,k\right\} $ such that $j\in \Gamma _i$ and $P\varphi
_{i...j}\circ P\pi _i\circ P\pi _{12...i}\circ P\pi
_{12...i+1}\circ ...\circ P\pi _{12...n-1}\left( \tau \right)
=P\varphi _{i...j}\circ P\pi _i\left( \tau \right) =\mu _j.$ By
the definition of the limit of a diagram we have that $P\varphi
_{i...j}\circ P\pi _i=P\pi _j$ for every $\left( i,j\right) \in
\Gamma $ and therefore $P\pi _j\left( \tau \right) =\mu _j.$ This
proves the lemma.
\end{proof}

{\bf 1)} First we restrict ourselves on the case that $X_i$ are finite. Since
$\chi$ is an affine surjective map (see Lemma \ref{l:2}) of the compact convex
polyhedron $P\left( {\lim }{D}\right) $ onto the compact convex polyhedron $
{\lim }P({D})$, the map $\chi$ is open (see \cite{Z}).

{\bf 2)} The case of zero-dimensional spaces $X_i$ is treated similarly as in
\cite{Z}. {\bf 3)} Suppose now that $X_i$ are arbitrary compact metrizable
spaces. There exist zero-dimensional compact spaces $X_i^{\prime }$ and
continuous surjective maps $f_i:X_i^{\prime }\rightarrow X_i$. Apply for the
diagram
\[
\xymatrix{P\left({\lim } {D'}\right)\ar[r]^{\chi^{\prime }
}\ar[d]_{P\left( h\circ \overset{k}{\underset{i=1}{\prod }}f_i\circ h^{-1}\right)
}&{\lim }P\left({D'}\right)\ar[d]^{h\circ \overset{k}{\underset{i=1}{\prod
}}P\left( f_i\right) \circ h^{-1}}\\ P\left( {\lim }{D}
\right)\ar[r]_\chi&{\lim }P({D})}
\]the lemma 2 and the lemma 2.1 from [1] together with the fact that
the characteristic map $\chi^{\prime }$ is open we see that the
map $\chi$ is open as well.
\end{proof}

Now we are going to prove the main result.

\begin{thm}
The probability measure functor $P$ is open-multicommutative.
\end{thm}

\begin{proof}
Consider an arbitrary cone $\left\{ T, h_i, i=1...n\right\}$ of the diagram $D$
such that the characteristic map
\[
\chi_{T,D}:T\rightarrow
{\lim}{D}
\] is surjective and open. It was proved before that the map
\[
\chi :P\left( {\lim }{D}
\right) \rightarrow {\lim }P\left( D\right)
\]is surjective and open. Since the functor $P$ preserves open and
surjective maps, the map
\[P\left( \chi_{T,D}
\right):P\left(T \right)\rightarrow P\left(
{\lim}{D}\right)
\]is also open and surjective. Consider the composition $\chi\circ P\left(
\chi_{T,D} \right)$. For every measure $\nu \in P\left({\lim }
{D} \right)$ we have $\chi\left(\nu\right)=\left(\nu_1,...,\nu_n\right)$ where
$\nu_i=P\pi_i\left(\nu\right) \in P\left( X_i\right)$. Let $\mu
\in P\left( T \right)$ and $\varphi \in C\left( T\right)$. It
holds
\begin{eqnarray*}
P\pi_i\circ P\left( \chi_{T,D}
\right)\left(\mu\right)\left(\varphi\right)&=&P\left( \chi_{T,D}
\right)\left(\mu\right)\left(\varphi\circ\pi_i\right)\\&=&\mu\left(\varphi\circ\pi_i\circ\chi_{T,D}\right)=
\mu\left(\varphi\circ
h_i\right)=Ph_i\left(\mu\right)\left(\varphi\right).
\end{eqnarray*}This implies that $\chi\circ P\left( \chi_{T,D}
\right)=\chi_{P\left(T\right),P\left(D\right)}$. Thus the
characteristic map $\chi_{P\left(T\right),P\left(D\right)}$ of the cone $\left\{
P\left(T\right), Ph_i, i=1...n\right\}$ is open and surjective. Since the cone is
arbitrarily chosen, this implies that the functor $P$ is open-multicommutative.

\end{proof}


\begin{thebibliography}{99}

\bibitem{Sh}  E.V. Shchepin, {\em Functors and uncountable powers of compacta.}
Uspekhi Mat. Nauk. 36, no.3, 3-62 (Russian)

\bibitem{DE} S. Z. Ditor, L. Q. Eifler, {\em Some open mapping theorems for measures},
Trans. Amer. Math. Soc. {\bf 164}(1972), 278--293.


\bibitem{b}
J. Bergin, {\em On the continuity of correspondences on sets of measures with
restricted marginals}, Economic Theory {\bf 13} (1999), 471-481.

\bibitem{Z}  M.M. Zarichnyi, {\em Correspondences of probability measures
with restricted marginals revisited}, Preprint.

\end{thebibliography}
\end{document}